\documentclass[titlepage,11pt]{article}
\usepackage{latexsym}
\usepackage{amsfonts}
\usepackage{amsmath}
\usepackage{epsfig}
\usepackage{pdfpages}
\usepackage{fancyhdr}
\newcommand\blackslug{\hbox{\hskip 1pt \vrule width 4pt height 8pt depth 1.5pt
        \hskip 1pt}}
\newcommand\bbox{\hfill \quad \blackslug \bigbreak}

\def\l{,\ldots,}

\title{Three-edge-colouring doublecross cubic graphs}
\author{Katherine Edwards\thanks{Supported by an NSERC PGS-D3 Fellowship and a Gordon Wu Fellowship.}\\
Princeton University, Princeton, NJ 08544
\\
\\
Daniel P. Sanders\thanks{Research performed while Sanders was a faculty member at Princeton University.}\\
Renaissance Technologies LLC,
East Setauket, NY 11733
\\
\\
Paul Seymour\thanks{Supported by ONR grants N00014-10-1-0680 and N00014-14-1-0084, and NSF
grant DMS-1265563.}\\
Princeton University, Princeton, NJ 08544
\\
\\
Robin Thomas\thanks{Supported by NSF grant number DMS-1202640.}\\
Georgia Institute of Technology, Atlanta, GA 30332}

\date{May 2, 2014; revised \today}

\newtheorem{thm}{}[section]

\newcommand{\Proof}{\noindent{\bf Proof.}\ \ }

\begin{document}
\maketitle
\begin{abstract}
A graph is {\em apex} if there is a vertex whose deletion makes the graph planar, and 
{\em doublecross} if it can be drawn in the plane with only two crossings,
both incident with the infinite region in the natural sense.
In 1966, Tutte~\cite{Tutte} conjectured that every two-edge-connected cubic graph with no
Petersen graph minor is three-edge-colourable. With Neil Robertson, two of us showed that
this is true in general if it is true for apex graphs and
doublecross graphs~\cite{RST1,RST2}. In another paper~\cite{apexcase}, two of us solved the apex case,
but the doublecross case remained open. Here we solve the doublecross case; that is, we prove that every two-edge-connected doublecross
cubic graph is three-edge-colourable.
The proof method is
a variant on the proof of the four-colour theorem given in~\cite{RSST}.

\end{abstract}

\section{Introduction}

In this paper, all graphs are finite and simple. The four-colour theorem (4CT),
that every planar graph is four-colourable, was proved by Appel and 
Haken~\cite{A&H1,A&H2} in 1977, and a simplified proof was given in 1997 by
three of us, with Neil Robertson~\cite{RSST}.  The 4CT is 
equivalent to the statement that every two-edge-connected cubic planar graph is
three-edge-colourable, and a strengthening was proposed as a conjecture 
by Tutte in 
1966~\cite{Tutte}; that every two-edge-connected cubic graph with no Petersen graph minor is three-edge-colourable. 
This paper is a step in the proof of Tutte's conjecture.

A graph $G$ is {\em apex} if $G\setminus v$ is planar for some vertex $v$ 
(we use $\setminus$ to denote deletion); and a graph $G$ is {\em doublecross}
if it can be drawn in the plane with only two crossings, both on the infinite 
region in the natural sense.

It is easy to check that apex and doublecross graphs do not contain the 
Petersen graph as a minor; but there is also a converse. Let us say a 
graph $G$ is {\em theta-connected}  
if $G$ is cubic and has girth at least five,
and $| \delta_{G} (X) | \geq 6$ for all
$X \subseteq V(G)$ with $|X| ,|V(G)  \setminus X| \geq 6$.
($ \delta_{G} (X) $ denotes the set of edges of $G$ with one end in $X$ and 
one end in $V(G)\setminus X$.)
Two of us (with Robertson) proved in~\cite{RST1} that every theta-connected graph with no Petersen graph minor is
either apex or doublecross (with one exception, that is three-edge-colourable);
and in~\cite{RST2} that every minimal 
counterexample to Tutte's conjecture was either apex or theta-connected.
It follows that every minimal counterexample to Tutte's conjecture is either
apex or doublecross, and so to prove the conjecture in general, it suffices
to prove it for apex graphs and for doublecross graphs. 
Two of us proved in~\cite{apexcase}
that every two-edge-connected apex cubic graph is three-edge-colourable,
so all that remains is the doublecross case, which is the objective of this paper. Our main theorem is:
\begin{thm}\label{mainthm}
Every two-edge-connected doublecross cubic graph is three-edge-colourable.
\end{thm}

The proof method is by modifying the proof of the 4CT given in~\cite{RSST}. 
Again we give a 
list of reducible configurations (the definition of ``reducible'' has to be 
modified to accommodate the two pairs of crossing edges), and a discharging 
procedure to prove that one of these configurations must appear in every 
minimal counterexample (and indeed in every non-apex theta-connected doublecross graph).
This will prove that there is no minimal counterexample, and so the theorem 
holds. 
Happily, the discharging rules given in~\cite{RSST} still work without
any modification.

\section{Crossings}

We are only concerned with graphs that can be drawn in the plane with two 
crossings, and one might think that these are not much different from planar 
graphs, and perhaps one could just {\em use} the 4CT rather 
than going to all the trouble of repeating
and modifying its proof. For graphs with one crossing this is true: here is
a pretty theorem of Jaeger~\cite{jaeger} (we include a proof because we think it is of interest):

\begin{thm}\label{onecrossing}
Let $G$ be a two-edge-connected cubic graph, drawn in the plane with one 
crossing. Then $G$ is three-edge-colourable.
\end{thm}
\Proof Let $e,f$ be the two edges that cross one another, and let
$e = z_1z_3$ and $f = z_2z_4$ say. Let $H$ be obtained from $G$
by deleting $e,f$, adding four new vertices $y_1\l y_4$, and edges $g_i=y_iz_i$
for $1\le i\le 4$. Thus every vertex of $H$ has degree three, except
for $y_1\l y_4$ which have degree one; and $H$ can be drawn in a closed
disc with $y_1\l y_4$ drawn in the boundary of the disc in order.
By adding four edges $y_1y_2,y_2y_3,y_3y_4$ and $y_4y_1$, we obtain
a two-edge-connected cubic planar graph, which therefore is 
three-edge-colourable, by the 4CT. Consequently $H$ is also 
three-edge-colourable; let $\phi:E(H)\rightarrow\{1,2,3\}$ be a three-edge-colouring of $H$. Since each colour 
appears at every vertex of $H$
different from $y_1\l y_4$, and there are an even number of such vertices,
it follows that each of the three colours
appears on an even number of $g_1\l g_4$. In particular, if $\phi(g_1) = \phi(g_3)$
then $\phi(g_2) = \phi(g_4)$, 
giving a three-edge-colouring of $G$ as required. We may assume then that
$\phi(g_1)\ne \phi(g_3)$, and similarly $\phi(g_2)\ne \phi(g_4)$. From
the symmetry we may therefore assume that $\phi(g_i) = 1$ for $i = 1,2$,
and $\phi(g_i) = 2$ for $i = 3,4$. Let $J$ be the subgraph of $H$ with vertex 
set $V(H)$ and edge set all edges $e$ of $H$ with $\phi(e)\in \{1,2\}$.
Every vertex of $H$ different from $y_1\l y_4$ therefore has degree two in $J$,
 and $y_1\l y_4$ have degree one; and so two components of $J$ are paths with
ends in $\{y_1\l y_4\}$, and all other components are cycles. Let the two 
components which are paths be $P_1,P_2$; and we may assume that $y_1$ is an end
of $P_1$. The second end of $P_1$ cannot be $y_3$, since then
$P_2$ would have ends $y_2,y_4$, which is impossible by planarity.
So the second end of $P_1$ is one of $y_2,y_4$; and in either case, if we 
exchange colours 1 and 2 on the edges of $P_1$, and otherwise leave $\phi$
unchanged, we obtain a new three-edge-colouring $\phi'$ of $H$, in which 
$\phi'(g_1) = \phi'(g_3)$ and 
$\phi'(g_2) = \phi'(g_4)$, 
which therefore gives a three-edge-colouring of $G$. This proves \ref{onecrossing}.~\bbox

\bigskip

We have tried (hard!) to do something similar to handle the doublecross case,
but failed; it seems necessary to do it the long way, modifying the proof of
the 4CT. Fortunately that is not as difficult as it was for the apex case in~\cite{apexcase}.

\section{A sketch of the proof}

Our proof is a modification of the proof of the 4CT, so let us begin with that. We can work in terms
of planar cubic graphs or in terms of planar triangulations, but we found it most convenient to 
use triangulations (because for instance, it is much
easier to present long lists of reducible configurations as subgraphs
of triangulations than as subgraphs of cubic graphs). 
Thus we need to show that every planar triangulation is 4-colourable. If there is a counterexample, 
the smallest one is
``internally 6-connected''; that is, every vertex-cutset of order at most five
separates at most one vertex from the rest of the graph. Now the average degree of 
vertices is just less than six, and so there are vertices of degree five; but we need to show that there is 
a small clump of vertices with small degrees. (We listed explicitly exactly what counts as a ``clump''; there were 633
of them in that paper. They have somewhere between 
four and twelve vertices all with small degree).
To show this we use a ``discharging procedure''; we initially assign ``charges'' (integers) to the vertices, depending 
on their degrees and with positive total; and then move small amounts among neighbouring vertices according
to some carefully chosen rules; and prove a theorem that any vertex that ends up with positive charge (and there
must be such a vertex) belongs to or neighbours one of the clumps we are hoping for.

Next we show that for any counterexample to the 4CT, there must be a smaller counterexample. We just showed that
it would contain one of the 633 clumps, and the clumps were carefully chosen with a special property:
for every clump $C$ in a internally 6-connected triangulation $T$, we can excise
$C$ from $T$ and replace it by something smaller, so that $T$ is changed into a smaller triangulation that is no easier to 4-colour.
This would show that the 4CT is true.

How can we show that for each of the 633 clumps, if it occurs in a triangulation, it can be replaced by something 
smaller without making the triangulation easier to colour?
The idea is, we show this clump by clump, using a computer. Let $C$ be one of our clumps, in an internally 6-connected 
triangulation $T$ that is not 4-colourable. (Since we are working through the clumps one at a time, at this stage 
we know $C$ completely, but have no knowledge of the remainder of $T$, which is some hypothetical triangulation.) 
Let $R$ be the cycle of $T$
(the ``ring'' of $C$) that bounds the region that we obtain when $C$ is deleted. 
(It might not actually be a cycle, but ignore that. We do know its length, because it is determined by the degrees in $T$ of
the  vertices in $V(C)$, and we know these numbers.)
Some 4-colourings of the ring can be extended to the part of $T$ outside
of $R$, and some extend inwards to $C$, but none of them extends in both directions because $T$ is not 4-colourable.
Since we know $C$, we can compute exactly which 4-colourings of $R$ extend inwards to $C$; so we know that the set
of ring-colourings that extend outwards is a subset of the complement of this. Also, the set that extend outwards 
has to be closed under
Kempe chain recolouring (we call this ``consistent''); and sometimes we find that this is impossible, that there is no such 
consistent set (except the null set,
which is consistent but cannot be the right answer). In this case, deleting $C$ still gives a planar graph that cannot be 4-coloured;
we don't need to substitute anything for it (unless we want to make a triangulation, when we had better put something in, just 
to make all the regions triangles). For other clumps, we find that there is a nonempty consistent set of 4-colourings of $R$
none of which extends inwards, and we figure out the maximal such set $S$, which is unique (fortunately); and 
$S$ is still very restricted, so much so that we can replace $C$ by some smaller gadget, with
the crucial property that no colourings in $S$ extend inwards to the gadget.

That is roughly the idea of the proof of the 4CT. How can we adapt it to handle doublecross graphs? 
The first issue is that we can't take the planar dual any more to make a triangulation. But we could if we modify the four 
crossing edges somehow to make the graph planar first; eg subdivide all four of them and identify their midpoints. So we have
something that is ``almost'' a planar dual of the cubic graph that we want to three-edge-colour. It is a planar triangulation
except for one region of length eight. This is worth doing, because
it is easier to work with triangulations than with cubic graphs.

Now we have $T$, which is almost a triangulation. If it really were a triangulation, it would have to contain
one of the 633 clumps that we used for the 4CT, but
\begin{itemize}
\item it isn't really a triangulation;
\item if it did contain one of those clumps, we need to make sure that the region of length eight is not involved;
\item even if we find a clump far away from the big region, it will not necessarily do what we want.
\end{itemize}
The third problem above arises because we have to work with Kempe chains in the original cubic graph, which is not planar;
now Kempe chains might
use the four crossing edges, and so we can't prove that the set of colourings of the ring that extend outwards is consistent 
in the old planar sense. All we know is, it must be consistent in a slightly weaker sense that allows for a little bit 
of crossing between the various Kempe chains. Call this new sense ``XX-consistent''. 
The maximal XX-consistent subset of ring-colourings disjoint from those that extend inwards gets bigger than in the planar case; 
and it might now contain
a colouring that extends inwards to the gadget we were planning to replace the clump by.  We need to search for
another set of clumps to replace the 633, that work with this weaker definition of consistent. We did this; so the third
problem goes away, if we replace the 633 with the 756 clumps described at the end of the paper. And one really nice feature;
the same discharging procedure as for the 4CT still proves that one of these clumps must always be present. But we still have 
the first and second problems to handle. 

Here is a way to handle them. At the moment our ``almost''-triangulation $T$ has one special region of length eight, inside which the
four crossing edges were drawn. Fill in this region with a dense triangulation, so that $T$ is extended to something that really
is a triangulation, $T'$; and do it in such a way that all the new vertices close to the eight boundary vertices have big degree 
(at least 12). We know that $T'$ must contain one of the 756 clumps; but we want to get one of them inside of $T$. The idea is,
the discharging argument tells us a little more; that {\em every} vertex that ends up with positive charge belongs to or neighbours
one of the 756 clumps.
So if only we can prove that some vertex of $T$ itself ends up with positive charge, we would have one of the 756 at least neighbouring
a vertex of $T$; and it cannot be part in $T$ and part sticking out, because all the new vertices close to the boundary have big 
degree. So it would be completely in $T$, which is what we want.

How can we show that some vertex in $T$ ends up with positive charge? We need to show that the total initial charge assigned to 
the interior of $T$ is big, and not too much gets shipped over the border out of $T$. This boils down to something easy. 
Take the initial doublecross cubic graph that we want to three-edge-colour, delete the four crossing edges, and look at the 
length of the cycle $Z$ that bounds the region in which the four crossing edges were drawn. Provided $Z$ has length at least
21, the discharging argument will give us what we want, one of the 756 clumps in the right place where it can be used.
(We are back-and-forth between the cubic graph and its ``dual'' here, unfortunately.)

What do we do if $Z$ has length less than 21? In this case we need a completely different approach; we look for a ``clump'' 
in the cubic graph (perhaps we should call this a dual clump)
containing the four crossing edges and some edges of $Z$. We prove that there must always be one that
can be replaced by something smaller without making the graph 3-edge-colourable. Again, it is an argument using consistency under
Kempe chains, but here consistency is back to its
old sense from the proof of the 4CT, because all the nonplanarity is inside the clump itself, and the Kempe chains are out 
in the planar part of the graph. And it works, just -- we could do this up to length 20, but not for 
length 21, which is exactly the cutoff for the discharging procedure to work.

\section{XX-good configurations and a discharging argument}

It is most convenient to present the reducible configurations as subgraphs
of triangulations, even though most of the argument is done in terms of cubic graphs.
A {\em drawing} is defined as in~\cite{RSST}, and therefore has no ``crossings''. 
(Sometimes we need to allow crossings, but then we say so.)
A {\em triangulation} $T$ means a non-null drawing in a 2-sphere such that every
region is bounded by a cycle of length three (briefly, ``is a triangle''), and a
{\em near-triangulation} is a non-null connected drawing in the plane
such that every finite region is a triangle.
If $T$ is a near-triangulation, its infinite region is bounded by a cycle if and only if $T$ is two-connected; and
if so, we denote this cycle by $T_\infty$.
A {\em configuration} $K$ consists of a near-triangulation
$G_K$ together with a map $\gamma_K:V(G_K)\rightarrow \mathbb{Z}$ 
($\mathbb{Z}$ denotes the set of all integers) with the following properties, where $d(v)$ is the degree of $v$ in $G_K$:
\begin{itemize}
\item $|V(G_K)|\ge 2$;
\item for every vertex $v$, $G_K\setminus v$ has at most two components, and if there are two, then
$\gamma_K(v) = d(v)+2$; 
\item for every vertex $v$, if $v$ is not incident with the infinite region then $\gamma_K(v) = d(v)$; and
otherwise $\gamma_K(v)>d(v)$, and in either case $\gamma_K(v)\ge 5$;
\item $K$ has ring-size at least two, where the {\em ring-size} of $K$ is defined to be
$\sum_v(\gamma_K(v)-d(v)-1)$, summed over all vertices $v$ incident with the infinite region such that
$G_K\setminus v$ is connected. (In fact, all configurations used in this paper have ring-size at least six.)
\end{itemize}
We use the same conventions as in~\cite{RSST} to describe configurations, and in particular, we use the
same vertex shapes in drawings to represent the numbers $\gamma_K(v)$. Thus, a vertex $v$ with $\gamma_K(v)=5$ is represented by
a solid circle; 7 by a hollow circle; 8 by a square; 9 by a triangle; 10 by a pentagon; and 6 and 11
by a point (there is only one instance of 11, in the last configuration of the list).
\begin{figure}[h]
\begin{center}
\includegraphics[scale=0.17]{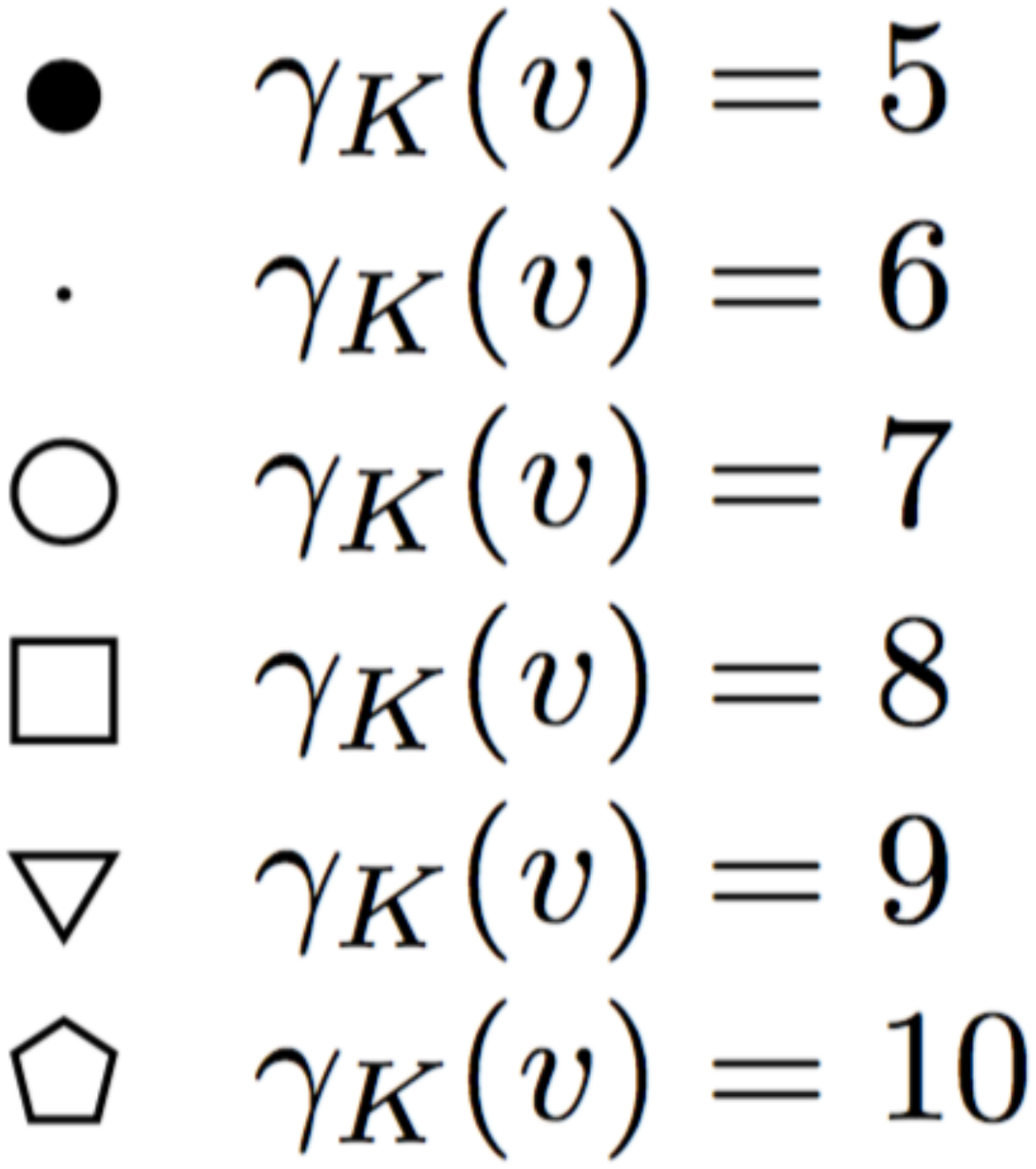}
\end{center}
\end{figure}
Two configurations $K,L$ are {\em isomorphic} if there is a homeomorphism of the plane mapping $G_K$ to $G_L$ and
mapping $\gamma_K$ to $\gamma_L$.
In the appendix to this paper
there are 756 configurations. Any configuration isomorphic to one of these is called an {\em XX-good configuration}. Note that
every XX-good configuration $K$ has the property
that $\gamma_K(v)\le 11$ for every $v\in V(G_K)$.

A triangulation $T$ (in a 2-sphere $\Sigma$) or a near-triangulation (in a plane $\Sigma$)
is {\em internally six-connected} if for every cycle $C$ of $T$ with length at most five,
either some open disc in $\Sigma$, bounded by $C$, contains no vertex of $T$, or $C$ has length five
and some such open disc contains a unique vertex of $T$.

We say a configuration $K$ {\em appears} in a triangulation $T$ if
\begin{itemize}
\item $G_K$ is an induced subgraph
of $T$; and
\item for each $v\in V(G_K)$, $\gamma_K(v)$ equals the degree of $v$ in $T$.
\end{itemize}

Let $T$ be a two-connected near-triangulation.
We say a configuration $K$ {\em appears internally} in $T$ if
\begin{itemize}
\item $G_K$ is an induced subgraph
of $T\setminus V(T_\infty)$;
\item for each $v\in V(G_K)$, $\gamma_K(v)$ equals the degree of $v$ in $T$; and
\item every vertex or edge of $T$ that does not belong to $G_K$ is drawn within the infinite region of $G_K$.
\end{itemize}

The main result of this section is the following.

\begin{thm}\label{unav}
Let $J$ be a two-connected, internally six-connected near-triangulation. Suppose that
$J_\infty$ is an induced subgraph of $J$, and
there are at least $4|V(J_\infty)|-11$ edges in $J$ between $V(J_\infty)$ and $V(J)\setminus V(J_\infty)$.
Then some XX-good configuration appears internally in $J$.
\end{thm}

To prove this, we use
the discharging procedure from~\cite{RSST}, so now we turn to that.
A {\em discharging function} in a
triangulation $T$ means a map $\phi$ from the set of all ordered pairs of adjacent vertices of $T$
into $\mathbb{Z}$, such that $\phi(u,v)+\phi(v,u) = 0$ for all adjacent $u,v$.
In~\cite{RSST}, we defined an explicit discharging function, for every internally six-connected triangulation.
Since it is rather complicated, we refer the reader to~\cite{RSST} for the details of its definition.
We need the following two properties; the first can be verified by hand, but the second needs a computer.

\begin{thm}\label{discharge}
Let $T$ be an internally six-connected triangulation, and let $\phi$ be the discharging function
defined in {\rm \cite{RSST}}. Then
\begin{itemize}
\item for every edge $uv$, if $\phi(u,v)>5$ then some XX-good configuration appears in $T$ and contains $u$
\item for every vertex $u$ of $G$, if $10(6-d_T(u))>\sum_v \phi(u,v)$ (where the sum is
over all vertices $v$ adjacent to $u$) then some XX-good configuration $K$ appears in $T$, and moreover
either $u$ or some neighbour of $u$ is a vertex of $G_K$.
\end{itemize}
\end{thm}

Both of these statements are minor variants of theorems proved in~\cite{RSST} 
(theorems 4.7 and 4.4 of that paper, respectively) and the 
methods of proof are unchanged.
The proof of the first statement is virtually identical with the proof of theorem 4.7 of~\cite{RSST}, because
all the ``good configurations'' used in that proof are also XX-good, except for one, the configuration
called conf(2,10,6) in that paper, which is not XX-good. This is only needed at one step of the proof, and 
at that step we can use a different configuration instead, which is XX-good, the second on line 5 of page 3 of the Appendix.
In the proof in~\cite{RSST} of theorem 4.7 of that paper, we did not include in the statement of the theorem
that the good configuration we find contains $u$,
as we are claiming now; but in fact the proof shows that.

The proof 
of the second needs analogues of theorems 4.5, 4.6, 4.8, 4.9 of the same paper. Again, in~\cite{RSST} we did not 
include in the statement of the theorem that the good configuration we find contains $u$ or one of its neighbours;
but this is implied by the fact that the configuration is always found within a ``cartwheel'' centred on $u$.
The proofs of the analogues of 4.5, 4.6 and 4.8 are unchanged, because all the good
configurations used for those proofs in~\cite{RSST} are also XX-good. The analogue
of 4.9 is proved by computer. The computer program just checks a machine-readable proof
of unavoidability, and is
the same as was used in~\cite{RSST}; we just changed its two inputs, the list
of configurations we want to prove unavoidable, and the files containing the machine-readable proofs. We are making
the program and the computer-readable proofs available on the
arXiv~\cite{programs}. 

\bigskip

\noindent{\bf Proof of \ref{unav}.\ \ }
Let $C=J_\infty$, and
take a drawing of $J$ in a 2-sphere such that $C$ bounds some region $r_0$. 
Thus one region of $J$ has a boundary of length $|V(C)|$, and all others have length three.
Let $J$ have $r$ regions; then by Euler's formula, $|V(J)|-|E(J)|+r = 2$, and so
$$\sum_{v\in V(J)}6-d_J(v) = 6|V(J)|-2|E(J)|=6(|E(J)|+2-r)-2|E(J)| = 4|E(J)|-6r+12.$$
The sum of the lengths of the regions of $J$ is
$|V(C)|+3(r-1)$, and so
$2|E(J)| = |V(C)|+3(r-1)$. We deduce that
$$\sum_{v\in V(J)}6-d_J(v) =4|E(J)|-6r+12= 2(|V(C)|+3(r-1))-6r+12= 2|V(C)|+6.$$

Let there be $k$ edges between $V(J)\setminus V(C)$ and $V(C)$.
It follows that $\sum_{v\in V(C)}d_J(v) = k+2|V(C)|$.
Consequently
$$\sum_{v\in V(C)} 6-d_J(v)= 6|V(C)|-(k+2|V(C)|)=4|V(C)|-k$$
and so
$$\sum_{u\in V(J)\setminus V(C)} 6-d_J(u)= (2|V(C)|+6)-(4|V(C)|-k) = k+6-2|V(C)|.$$

Now extend the drawing of $J$
to a drawing of an internally six-connected triangulation $T$, by adding more vertices and edges
drawn within $r_0$, in such a way that every vertex in $V(C)$ has degree in $T$ at least 12. (It is easy
to see that this is possible, since $J$ is internally six-connected.)
\\
\\
(1) {\em If there is an XX-good configuration $K$ that appears in $T$ such that some vertex in $V(J)\setminus V(C)$ either
belongs to $G_K$ or has a neighbour in $G_K$, then $K$ appears internally in $J$.}
\\
\\
No vertex in $V(C)$ belongs to $V(G_K)$, since $\gamma_K(v)\le 11$ for every vertex $v$ of $K$, and
$d_T(v)\ge 12$ for every $v\in V(C)$. From the hypothesis, it follows that some vertex of $G_K$ belongs to $V(J)\setminus V(C)$,
and hence $V(G_K)\subseteq V(J)\setminus V(C)$ since $G_K$ is connected. Moreover, every finite region of $G_K$ is a finite region
of $J$, since $J$ is internally six-connected. But every vertex in $V(J)\setminus V(C)$ has the same degrees 
in $T$ and in $J$,
and so $K$ appears internally in $J$. This proves (1).

\bigskip

Let $\phi$ be
the discharging function on $T$ defined in~\cite{RSST}.
Suppose first that $\phi(u,v)>5$ for some edge $uv$ of $T$ with $u\in V(J)\setminus V(C)$
and $v\in V(C)$. By the first statement of \ref{discharge}, some XX-good configuration
$K$ appears in $T$ and contains $u$, and the result follows from (1).

Thus we may assume that there is no such edge $uv$. Let there be $k$ edges $uv$ in $J$ with $u\in V(J)\setminus V(C)$
and $v\in V(C)$.
Consequently the sum of $\phi(u,v)$,
over all edges $uv$ of $T$ with $u\in V(J)\setminus V(C)$
and $v\in V(C)$, is at most $5k$.
But this equals the sum over all $u\in V(J)\setminus V(C)$, of the sum of $\phi(u,v)$
over all neighbours $v$ of $u$, since $\phi(u,v) = -\phi(v,u)$ for all $u,v$. Therefore
the sum over all $u\in V(J)\setminus V(C)$ of
$$10(6-d_J(u))-\sum_{uv\in E(J)} \phi(u,v)$$
is at least $10(k+6-2|V(C)|)-5k$.
Since $k\ge 4|V(C)|-11$ by hypothesis, the last is positive, and so
there exists $u\in  V(J)\setminus V(C)$ such that
$$10(6-d_J(u))-\sum_{uv\in E(J)} \phi(u,v)>0.$$

For such a vertex $u$, its degrees in $J$ and $T$ are the same; and so by the second assertion
of \ref{discharge}, there is an XX-good configuration $K$ that appears in $T$, such that
either $u$ or some neighbour of $u$ is a vertex of $G_K$.
But then again, the result follows from (1).
This proves \ref{unav}.~\bbox

\section{The doublecross edges}

By a {\em minimal counterexample}, we mean a cubic two-edge-connected doublecross graph $G$ which is
not three-edge-colourable, and such that 
every cubic two-edge-connected minor of $G$ is three-edge-colourable except $G$ itself.

\begin{thm}\label{properties}
Let $G$ be a minimal counterexample. Then 
\begin{itemize}
\item $G$ is theta-connected;
\item 
there are four edges $g_1\l g_4$ of $G$,
and the graph $G\setminus \{g_1,g_2,g_3,g_4\}$ can be drawn in the plane such that its infinite region
is bounded by a cycle $Z$;
\item there are eight vertices $z_1\l z_8$ of $Z$, distinct and in order, such that $g_1=z_1z_3$,
$g_2 = z_2z_4$, $g_3 = z_5z_7$ and $g_4 = z_6z_8$.
\end{itemize}
\end{thm}
\Proof
By the result of~\cite{apexcase}, $G$ is not apex since it is not three-edge-colourable.
Since $G$ is also a minimal counterexample to Tutte's conjecture, the result of~\cite{RST2}
implies that $G$ is theta-connected. Now
$G$ can be drawn in the plane with only two crossings both on the infinite region.
Let the crossing pairs of edges be $(g_1,g_2)$ and $(g_3,g_4)$. Since $G$ is not apex, it follows
that 
$g_1,g_2\ne g_3,g_4$, and indeed $g_1,g_2$ are disjoint from $g_3,g_4$.
If $g_1$ shares an end with $g_2$,
then the drawings of $g_1,g_2$ can be rearranged to eliminate their crossing,
and again $G$ is apex, which is again impossible (indeed, in this case
$G$ has crossing number at most one, and so we could use \ref{onecrossing} instead
of the result of~\cite{apexcase}). So $g_1,g_2,g_3,g_4$ are disjoint. 
The graph obtained
by deleting the edges $g_1\l g_4$ is two-edge-connected since $G$ is theta-connected, and it is drawn in
the plane such that $g_1,g_2,g_3,g_4$ are all drawn within its infinite region. Since
it is two-edge-connected, its infinite region is bounded by a cycle.
This proves \ref{properties}.~\bbox

\bigskip

To prove \ref{mainthm}, we use the same approach as the proof of the four-colour theorem, proving the existence
of an unavoidable set of reducible subgraphs. Some of these reducible subgraphs use all four of the edges
$g_1\l g_4$, and the others use none of them. The length of the cycle $Z$ of \ref{properties} is the deciding factor
here; if $|E(Z)|\le 20$ we will show the existence of a reducible subgraph using $g_1\l g_4$, while if $|E(Z)|\ge 21$
we will show the presence of one of the other kind. In this section we handle the case when $|E(Z)|\le 20$.

\begin{thm}\label{longcycle}
If $G$ is a minimal counterexample and $Z$ is the cycle as in \ref{properties}, then $|E(Z)|\ge 21$.
\end{thm}
\Proof 
Let $z_1\l z_8$ be as in \ref{properties}, and for $1\le i\le 8$ let $Z_i$ be the path of $Z$
with ends $z_i,z_{i+1}$ that contains no other member of $\{z_1\l z_8\}$ (where $z_9$ means $z_1$).
For $1\le i\le 8$, let $L_i$ denote $|E(Z_i)|$. We observe:
\begin{itemize}
\item $L_1\l L_8\ge 1$, because $z_1\l z_8$ are all distinct;
\item $L_1+L_2\ge 4$ since every cycle of $G$ has length at least five; and for the same reason
$L_2+L_3, L_5+L_6, L_6+L_7\ge 4$ and $L_1+L_3, L_5+L_7\ge 3$;
\item $L_1+L_2+L_3\ge 7$ since there are at least six edges with exactly one end
in $V(Z_1\cup Z_2\cup Z_3)$ (because $G$ is theta-connected); and similarly $L_5+L_6+L_7\ge 7$;
\item $L_1+L_8\ge 3$, because if $Z_1,Z_8$ both have length one then $G$ is apex (deleting the end of $Z_8$ not in $Z_1$
makes the graph planar); and similarly $L_3+L_4, L_4+L_5, L_7+L_8\ge 3$.
\end{itemize}

A choice of the 8-tuple $(L_1\l L_8)$ is called {\em plausible} if it satisfies the conditions just listed.
Suppose that $|E(Z)|\le 20$; then there are only finitely many plausible choices for $(L_1\l L_8)$,
and we handle them one at a time. 
Now, therefore, we assume that we are dealing with some such plausible choice, and so we know the lengths $L_1\l L_8$.
Let $G^-$ be the graph obtained from $G$ by deleting the four crossing edges
$g_1 = z_1z_3$, $g_2=z_2z_4$, $g_3=z_5z_7$ and $g_4=z_6z_8$. Then $Z$ is a cycle of $G^-$, bounding a region in a planar 
drawing of $G^-$. Every vertex of $Z$ different from $z_1\l z_8$ is incident with an edge of $G^-$ that
does not belong to $E(Z)$. Let the vertices of $Z$ different from $z_1\l z_8$ be $v_1\l v_k$ say,
numbered in circular order (starting from some arbitrary first vertex), and for $1\le i\le k$
let $f_i$ be the edge of $G$ incident with $v_i$ and not in $E(Z)$. Note that $f_1\l f_k$
might not all be distinct (because for instance some $f_i$ might be incident with a vertex of
the interior of $Z_2$ and incident with a vertex of the interior of $Z_6$). Let $F=\{f_1\l f_k\}$,
and let $\mathcal{D}$ be the set
of all maps from $F$ to $\{1,2,3\}$.

A subset $\mathcal{C}\subseteq \mathcal{D}$ is said to be
{\em consistent} if it has the following property.
For all distinct $x,y\in \{1,2,3\}$, and each $\phi\in \mathcal{C}$, let $F_{x,y}$ be the 
set of all $f\in F$ with $\phi(f)\in \{x,y\}$; then there is a partition $\Pi$ of $F_{x,y}$ into sets
of size one and two, with the following properties:
\begin{itemize}
\item for $f\in F_{x,y}$, the member of $\Pi$ containing $f$ has size one if and only if both ends of
$f$ belong to $V(Z)$;
\item for $1\le a<b<c<d\le  k$, not both $\{f_a,f_c\}, \{f_b,f_d\}\in \Pi$; and
\item $\phi'\in \mathcal{C}$ for every subset $F'\subseteq F_{x,y}$ which is expressible as a union of members of $\Pi$,
where $\phi'$ is defined by
$$\phi'(f) = \begin{cases} \phi(f)&\mbox{if } f\in F\setminus F'\\
			y&\mbox{if } f\in F'$ and $\phi(f) = x\\
			x&\mbox{if } f\in F'$ and $\phi(f) = y.
	\end{cases}$$
\end{itemize}

For any graph $H$ with $F\subseteq E(H)$, we denote by $\mathcal{C}_H$ the set of all members of
$\mathcal{D}$ that can be extended to a three-edge-colouring of $H$.
Let $J= G^-\setminus E(Z)$, and let $K$ be the subgraph of $G$ formed by the edges in
$F\cup E(Z)\cup\{g_1,g_2,g_3,g_4\}$ and their incident vertices.
Since $|E(Z)|\le 20$, it follows that $|F|\le 12$.

For each plausible choice of $(L_1\l L_8)$ (except one, that we handle separately), 
there are three steps to be carried out 
on a computer, which we explain now. All three involve computation with subsets of $\mathcal{D}$,
but since $|F|\le 12$ all three steps are easily implemented. 
\begin{description}
\item {\bf Step 1:}  Compute $\mathcal{C}_K$.
\item {\bf Step 2:} Compute the maximal consistent subset $\mathcal{C}$ of
$\mathcal{D}\setminus \mathcal{C}_K$. (The union of any two consistent sets
is consistent, and so there is a unique maximal consistent subset of any set.)
\item {\bf Step 3:} Verify that 
there is a graph $K'$, obtained from $K$ by deleting one or two of the edges $g_1\l g_4$,
and suppressing the vertices of degree two that arise, such that
$\mathcal{C}\cap \mathcal{C}_{K'}=\emptyset$. (This involves looking at the handful of possibilities
for $K'$ and computing $ \mathcal{C}_{K'}$ for each of them.)
\end{description}

The exception is the 8-tuple $(4,1,4,1,4,1,4,1)$, which is plausible but for which there is no $K'$ as in step 3 above.
For this 8-tuple let $K'$ be obtained from $K$ by deleting the unique edge in $Z_2$ and suppressing the vertices of degree
two that arise. Again we check that $\mathcal{C}\cap \mathcal{C}_{K'}=\emptyset$.

In all cases, 
since $G$ is not three-edge-colourable, and $J\cup K = G$, it follows
that $\mathcal{C}_J$ and $\mathcal{C}_K$ are disjoint.
The planarity of $J$ implies that $\mathcal{C}_J$ is consistent (this is easy to see, and is a standard 
argument in proving the four-colour theorem -- see~\cite{RSST}). 
Consequently $\mathcal{C}_J\subseteq \mathcal{C}$, and so $\mathcal{C}_J\cap \mathcal{C}_{K'}=\emptyset$.
It follows that $J\cup K'$ is not three-edge-colourable. But $J\cup K'$ can be obtained from $G$
by deleting either one or two disjoint edges, and suppressing the vertices of degree two that arise;
and so $J\cup K'$ is two-edge-connected (since $G$ is theta-connected), doublecross, and smaller than our supposedly
minimal counterexample, which is impossible. 
This proves \ref{longcycle}.~\bbox

We have omitted the details of the computer checking; this is all straightforward. (The program is available on the arXiv.)
There are 2957 plausible
8-tuples to check, up to symmetry, but the program only takes about a minute to do them all, 
so we were content with that. If desired, running through 
all possible choices of $(L_1\l L_8)$ could be made more
efficient at the cost of complicating the proof. For instance, we could quickly dispose of the case when
$\min(L_1,L_2,L_3)=\min(L_5,L_6,L_7)=1$, because in this case $G$ contains a ``C-reducible'' subgraph (of a different kind),
no matter what the other six lengths are. But we are aiming for simplicity 
rather than speed here.

\section{Islands}

An {\em island} means a graph $I$ drawn in the plane, with the following properties:
\begin{itemize}
\item $I$ is two-connected;
\item every vertex has degree two or three; and
\item every vertex of degree two is incident with the infinite region.
\end{itemize}

Let $I$ be an island, and $J$ be a geometric dual, where $j\in V(J)$ corresponds to the infinite region of $I$.
For each $v\in V(J)\setminus \{j\}$, let $\gamma(v)$ be the length of the region of $I$ that corresponds to $v$.
Then the pair $(J\setminus \{j\}, \gamma)$ might or might not be a configuration; but more important,
for every configuration $K$, there is an island that gives rise to it in this way, unique up to homeomorphism
of the plane. (We leave checking this to the reader. One way is to go to the ``free completion'' of $K$ defined in
\cite{RSST}, take a dual, and delete the vertex corresponding to the infinite region.) 
We call this the {\em island of $K$}, and denote it by $I(K)$. 

We need to work mostly with the islands of the XX-good configurations, but it is more compact to draw the
configurations themselves. Sometimes we need to refer to an edge $e$ of one of the islands, say of $I(K)$.
Now $e$ corresponds to some edge $f$ of $J$ under the duality (where $J$ is as before), 
and if $f\in E(G_K)$ then we can refer to $e$
by defining it as the edge dual to $f$. But sometimes the edge $f$
is not an edge of $G_K$.
For this reason, in the list of XX-good
configurations, some vertices are drawn with extra ``half-edges''. These indicate
some of the edges of $J$ that are not edges of $G_K$, for convenience in referring to certain edges of $I(K)$.

\begin{thm}\label{unavisland}
Let $G$ be a minimal counterexample, and let $Z, g_1\l g_4$ be as in \ref{properties}.
Then there is a cycle of $G\setminus \{g_1\l g_4\}$, bounding a closed disc $\Delta$,
such that the subgraph of $G$ formed by the vertices
and edges drawn in $\Delta$ is an island of some XX-good configuration.
\end{thm}
\Proof
Let $G^-=G\setminus \{g_1\l g_4\}$; and let us extend the drawing of $G^-$ by adding one new vertex $z_\infty$ 
(drawn within the infinite region of $G^-$)
and eight new edges, joining $z_\infty$ to the eight ends of $g_1\l g_4$, forming $G^+$ say.
Thus in $G^+$, every vertex
has degree three except for $z_\infty$, which has degree eight. (We can think of $G^+$ as obtained from $G$
by subdividing the four edges $g_1\l g_4$ and identifying the four new vertices.)

Now take a geometric dual $T$ of $G^+$, such that  $z_\infty$ belongs to the infinite region of $T$. 
It follows that:
\begin{itemize}
\item $T$ is a near-triangulation;
\item the cycle $T_\infty$ is an induced subgraph of $T$;
\item $T$ is internally six-connected (since $G$ is theta-connected); and
\item the number of edges of $T$ between $V(T_\infty)$ and $V(T)\setminus V(T_\infty)$ is at least 21, since $|E(Z)|\ge 21$ 
by \ref{longcycle}, and
for each $e\in E(Z)$, the corresponding edge of $T$ has one end in $V(T_\infty)$ and the other in $V(T)\setminus V(T_\infty)$.
\end{itemize}

Since $|V(T_\infty)|=8$ and so
$ 4|V(T_\infty)|-11=21$,  \ref{unav} implies that some XX-good configuration $K$ appears internally in $T$.
Consequently the union of the closures of the regions of $G^-$ that correspond to vertices in $G_K$ is a closed disc that
defines an island satisfying the theorem. This proves \ref{unavisland}.~\bbox

\section{Reducibility}

It remains to show that the outcome of \ref{unavisland} is impossible, but for that we need to discuss
reducibility further.

If $a,b,c,d$ are integers and $1\le a<b<c<d$, we call $\{\{a,c\},\{b,d\}\}$ a 
{\em cross}.
Let $\Pi$ be a finite set of finite sets of positive integers, each of cardinality two and pairwise disjoint.  
We say that $\Pi$ is {\em doublecross} if the following conditions hold:
\begin{itemize}
\item at most two subsets of $\Pi$ are crosses 
\item if $A,B,C,D\in \Pi$ are distinct, and $\{A,B\},\{C,D\}$ are crosses, 
let $X=A\cup B\cup C\cup D$ (so $|X|=8$). Then for all 
$P\in \Pi$ with $P\cap X=\emptyset$  there do not exist $x_1,x_2\in X$ 
such that $\{P,\{x_1,x_2\}\}$ 
is a cross.
\end{itemize}
This is equivalent to the following geometric condition, which may be easier to grasp: choose $k\ge 3$ such that 
$A\subseteq \{1\l k\}$ for each $A\in \Pi$. Take a regular $k$-vertex polygon in the plane, with vertices $v_1\l v_k$
in order.
For each $A\in \Pi$ draw a line segment $L_A$
between $v_i,v_j$, where $A=\{i,j\}$.
Then we ask that
\begin{itemize}
\item no point of the plane belongs to more than two of $L_A\;(A\in \Pi)$;
\item at most two points of the plane
belong to more than one of these lines; 
\item if there are two points $x,y$ each belonging to two of the lines $L_A\;(A\in \Pi)$ say, 
then either some $L_A$ contains them both, or no $L_A$
intersects the interior of the line segment between $x,y$.
\end{itemize}
We leave the equivalence to the reader.

Let $k\ge 1$, and 
let $\mathcal{D}$ be the set of 
all maps $\phi:\{1\l k\}\rightarrow\{1,2,3\}$.
We say a subset $\mathcal{C}$ of $\mathcal{D}$ is {\em XX-consistent} if it has the following property.
For all distinct $x,y\in \{1,2,3\}$, and each $\phi\in \mathcal{C}$, let $R_{x,y}$ be the
set of all $i\in \{1\l k\}$ with $\phi(i)\in \{x,y\}$; then there is a doublecross partition $\Pi$ of 
$R_{x,y}$, such that
$\phi'\in \mathcal{C}$ for every subset $F'\subseteq R_{x,y}$ which is expressible as a union of members of $\Pi$,
where $\phi'$ is defined by
$$\phi'(f) = \begin{cases} \phi(f)&\mbox{if } f\in \{1\l k\}\setminus F'\\
                        y&\mbox{if } f\in F'$ and $\phi(f) = x\\
                        x&\mbox{if } f\in F'$ and $\phi(f) = y.
        \end{cases}$$

Let $I$ be the island of a configuration $K$. 
We say $F\subseteq E(I)$ is a {\em matching} if no two edges in $F$ have a common end, and $V(F)$ denotes the set of vertices
incident with an edge in $F$.
A {\em three-edge-colouring modulo $F$} of $I$ means
a map $\phi:\;E(I)\setminus F\rightarrow\{1,2,3\}$, such that for all distinct edges $e,f\in E(I)\setminus F$ 
with a common end $v$ say, $\phi(e)=\phi(f)$ if and only if $v\in V(F)$. 

Let the vertices with degree two in $I$ be $v_1\l v_k$ in order on the boundary of the infinite region.
With $\mathcal{D}$ as before, let $\mathcal{C}_K$ be the 
set of all $\psi \in \mathcal{D}$ such that there is a three-edge-colouring $\phi$ of $I$ with $\phi(e)\ne\psi(i)$
for $1\le i\le k$ and for each edge $e$ of $I$ incident with $v_i$.
We say that $K$ is {\em XXD-reducible} if there is no non-null
XX-consistent subset of $\mathcal{D}\setminus \mathcal{C}_K$. We say that $K$ is {\em XXC-reducible}
if there is a matching $F$ of $I$ with the following properties:
\begin{itemize}
\item $1\le |F|\le 4$.
\item If $|F|=4$, then either some finite region of $I$ is incident with at least three members of $F$,
or there are two finite regions of $I$, say $r,s$, such that some edge of $I$ is incident with both $r,s$,
and every edge of $F$ is incident with one of them.
\item Let $\mathcal{C}_F$ be the set of all $\psi \in \mathcal{D}$ such that there is a three-edge-colouring
modulo $F$ of $I$, say $\phi$, 
such that for $1\le i\le k$ and every edge $e\in E(I)\setminus F$, $\phi(e)= \psi(i)$ if and only if $v_i\in V(F)$.
Then every XX-consistent subset of 
$\mathcal{D}\setminus \mathcal{C}_K$ is disjoint from $\mathcal{C}_F$.
\end{itemize}
We call such a set $F$ a {\em reducer} for $K$. (We will show that if $K$ appears in a minimal counterexample $G$, then 
deleting from $G$ the edges in $F$ and suppressing the resultant vertices of degree two will make a smaller counterexample,
which is impossible; and so $K$ cannot appear.)
We need:
\begin{thm}\label{Cred}
For every XX-good configuration $K$, either $K$ is XXD-reducible, or $K$ is XXC-reducible.
\end{thm}
\Proof The proof uses a computer. Each XX-good configuration $K$ is drawn in the Appendix,
and in that drawing sometimes some edges are drawn thickened. If no edges are thickened
then we claim $K$ is XXD-reducible, and otherwise we claim it is XXC-reducible, and the corresponding 
reducer $F$ corresponds to the set of thickened edges and half-edges of $G_K$ under planar duality. 
To show this, for each XX-good configuration $K$ in turn, 
we carry out two steps:
\begin{description}
\item {\bf Step 1:}  Compute $\mathcal{C}_K$.
\item {\bf Step 2:} Compute the maximal XX-consistent subset $\mathcal{C}$ of
$\mathcal{D}\setminus \mathcal{C}_K$. (The union of any two XX-consistent sets
is XX-consistent, and so there is a unique maximal XX-consistent subset of any set.)
\end{description}
If  $\mathcal{C}$ is empty we have verified that $K$ is XXD-reducible and we stop here.
Otherwise, we carry out:
\begin{description}
\item {\bf Step 3:} Let $F$ be the set of edges of $I(K)$ that correspond under
geometric duality to the thickened edges and half-edges of $G_K$ given in the Appendix, and 
verify that $F$ is a reducer for $K$.
\end{description}

This is just the same process as in the proof of the four-colour theorem~\cite{RSST}, and is carried
out on a computer the same way; we omit further details. (Again, we are making
the program available on the arXiv~\cite{programs}.)

\section{Assembling the pieces}

Now we combine these various lemmas to prove \ref{mainthm}. 

\bigskip

\noindent{\bf Proof of \ref{mainthm}.\ \ }
Suppose the result is false; then there is a minimal counterexample $G$. Let $Z, g_1\l g_4$
be as in \ref{properties}, and let $G^-=G\setminus\{g_1,g_2,g_3,g_4\}$.
By \ref{unavisland}, there is a cycle $D$ of $G^-$, bounding a closed disc $\Delta$,
such that the subgraph of $G$ formed by the vertices
and edges drawn in $\Delta$ is an island $I$ of some XX-good configuration $K$ say.
Let $v_1\l v_k$ be the vertices of $D$ that have degree two in $I$, numbered in order on $D$,
and for $1\le i\le k$, let $e_i$ be the edge of $G$ incident with $v_i$ and not in $E(I)$.
Note that $e_1\l e_k$ need not all be distinct, because some $e_i$ might have both ends in $V(D)$.

Now let $\mathcal{D}$ be the set of all maps from $\{1\l k\}$ to $\{1,2,3\}$.
We say a map $\phi:E(G)\setminus E(I)\rightarrow\{1,2,3\}$ is a {\em three-edge-colouring of $\overline{I}$} if $\phi(e)\ne \phi(f)$
for every two distinct edges $e,f\in E(G)\setminus E(I)$ with a common end in $G$.
Let $\mathcal{C}$ be the set of $\psi\in \mathcal{D}$ such that there is a three-edge-colouring $\phi$ of $\overline{I}$
with $\psi(i) = \phi(e_i)$ for $1\le i\le k$.
We claim:
\\
\\
(1) {\em $\mathcal{C}$ is non-null and XX-consistent.}
\\
\\
Clearly it is non-null, from the minimality of $G$.
Let $\psi\in \mathcal{C}$, and choose a three-edge-colouring $\phi$ of $\overline{I}$
with $\psi(i) = \phi(e_i)$ for $1\le i\le k$.
Let $x,y\in \{1,2,3\}$ be different, and let $H$ be the subgraph of $G$ formed by the edges $e\in E(G)\setminus E(I)$ with $\phi(e)\in\{x,y\}$
and their ends. It follows that every component of $H$ is either a cycle, or a path with distinct ends both in $V(D)$.
Let $\Pi$ be the set of all $\{i,j\}$ such that $1\le i<j\le k$ and some component of $H$ is a path with end-edges $e_i,e_j$.
(Possibly $e_i=e_j$, and this path has only one edge.) Then we can switch colours $x,y$ on any subset of $\{1\l k\}$
that is expressible as a union of members of $\Pi$, by exchanging the colours $x,y$ on the corresponding components of $H$.
It remains to show that $\Pi$ is doublecross. To see this, note first that if $\{\{a,b\},\{c,d\}\}\subseteq \Pi$ is a cross,
and $P,Q$ are the components of $H$ with end-edges $e_a,e_b$ and $e_c,e_d$ respectively, then either
$P$ contains one of $g_1,g_2$ and $Q$ contains the other, or $P$ contains one of $g_3,g_4$ and $Q$ contains the other.
Since $G$ can be drawn with no crossing pairs of edges except $g_1,g_2$ and $g_3,g_4$, and the two crossings they
form are on a common region, it follows that $\Pi$ is doublecross. This proves (1).

\bigskip

Now let $\mathcal{C}_K$ be as in \ref{Cred}; that is, 
$\mathcal{C}_K$ is the
set of all $\psi \in \mathcal{D}$ such that there is a three-edge-colouring $\phi$ of $I$ with $\phi(e)\ne\psi(i)$
for $1\le i\le k$ and for each edge $e$ of $I$ incident with $v_i$.
\\
\\
(2) {\em $\mathcal{C}_K\cap \mathcal{C}=\emptyset$.}
\\
\\
For suppose that $\psi\in \mathcal{C}_K\cap \mathcal{C}$. Choose a three-edge-colouring $\phi_1$ of $\overline{I}$ such that
 $\psi(i) = \phi_1(e_i)$ for $1\le i\le k$. Choose a three-edge-colouring $\phi_2$ of $I$ such that $\phi_2(e)\ne\psi(i)$
for $1\le i\le k$ and each edge $e$ of $I$ incident with $v_i$. For each edge $e$ of $G$, let
$$\phi(e)=\begin{cases}\phi_1(e)&\mbox{if } e\notin E(I)\\
                        \phi_2(e)&\mbox{if } e\in E(I).
        \end{cases}$$

We claim that $\phi$ is a three-edge-colouring of $G$. For let $e,f\in E(G)$ be distinct with a common end $v$ say.
If $e,f\in E(I)$ then 
$$\phi(e)=\phi_2(e)\ne \phi_2(f)=\phi(f)$$
since $\phi_2$ is a three-edge-colouring of $I$; and similarly $\phi(i)\ne \phi(j)$ if $i,j\notin E(I)$. We
may therefore assume that $e\in E(I)$ and $f\notin E(I)$; and consequently $v$ is one of $v_1\l v_k$, say $v_i$.
From the choice of $\phi_1$ it follows that $\psi(i) = \phi_1(e_i)=\phi(f)$; and from the choice of $\phi_2$,
$\phi(e)=\phi_2(e)\ne\psi(i)$. It follows that $\phi(e)\ne \phi(f)$. This proves (2).

\bigskip

By (1) and (2), it follows that $K$ is not XXD-reducible. Since $K$ is XX-good, it is therefore XXC-reducible by \ref{Cred}.
Let $F\subseteq E(I)$ be a reducer.
Let $\mathcal{C}_F$ be the set of all $\psi \in \mathcal{D}$ such that there is a three-edge-colouring
modulo $F$ of $I$, say $\phi$,
such that for $1\le i\le k$ and every edge $e\in E(I)\setminus F$ incident with $v_i$, $\phi(e)= \psi(i)$ if and only if $v_i\in V(F)$.
From the definition of a reducer, it follows that
every XX-consistent subset of
$\mathcal{D}\setminus \mathcal{C}_K$ is disjoint from $\mathcal{C}_F$, and so
in particular, $\mathcal{C}\cap \mathcal{C}_F=\emptyset$, by (1) and (2).
\\
\\
(3) {\em The graph $G\setminus F$ has a cutedge.}
\\
\\
Suppose it does not. Then from the minimality of $G$, there is a map $\phi:E(G)\setminus F\rightarrow\{1,2,3\}$, such that
for all distinct edges $e,f\in E(G)\setminus F$ with a common end $v$, $\phi(e)= \phi(f)$ if and only if $v\in V(F)$.
(To see this, suppress the vertices of degree two.) For $1\le i\le k$, let $\psi(i)=\phi(e_i)$.
Then $\psi \in \mathcal{C}\cap \mathcal{C}_F$, which is impossible. 
This proves (3).
\\
\\
(4) {\em $|F|=4$, and there is a cycle $W$ of $G$ of length five, such that 
$F\subseteq \delta_G(V(W))$.}
\\
\\
Let $f_0$ be a cutedge of $G\setminus F$. Consequently there exists $Y\subseteq V(G)$,                           
such that $f_0\in \delta_G(Y)\subseteq F\cup \{f_0\}$.
By
replacing $Y$ by its complement if necessary, we may assume that $|Y|\le |V(G)\setminus Y|$.
Suppose first that $| \delta_G(Y)|\le 4$. Since $G$ is theta-connected, it follows that $|Y|\le 2$.
Since $F\cap \delta_G(Y)$ is a matching, and $G$ is three-connected, this
is impossible.
Thus  $| \delta_G(Y)|\ge 5$, and so $|F|=4$, and $\delta_G(Y)=F\cup \{f_0\}$.

Since $G$ is theta-connected and $|Y|\le |V(G)\setminus Y|$, it follows that $|Y|\le 5$.
Now $|Y|\ge 4$ since
$F$ is a matching; and $|Y|$ is odd since $|\delta_G(Y)|$ is odd; so $|Y|=5$. Since $|\delta_G(Y)|=5$ and so there are
five edges of $G$ with both ends in $Y$, it follows that there is a cycle $W$ of $G$ with $V(W)\subseteq Y$; 
and since $G$ is theta-connected, $W$ has length
five. This proves (4).

\bigskip
We remark that
the edges in $F$ all belong to $E(I)$, but some of the other six edges of $G$ 
with an end in the cycle $W$ of (4) might not belong to
$E(I)$, and indeed might not belong to $E(G^-)$. We recall that from the choice of $F$, we have:
\\
\\
(5) {\em Either there exists a finite region of $I$ incident with three edges in $F$, or there are two finite regions
$r,r'$ of $I$, such that some edge of $I$ is incident with both $r,r'$, and every edge in $F$ is incident with one of $r,r'$.}

\bigskip

Let $W$ be as in (4), with vertices $w_1\l w_5$ in order, and for $1\le i\le 5$ let $h_i$ be the edge of $G$ incident with
$w_i$ and not in $E(W)$, where $F=\{h_1\l h_4\}$.
\\
\\
(6) {\em $W$ contains at least one of $g_1\l g_4$.}
\\
\\
For suppose not; then $W$ is a cycle of $G^-$, and consequently
bounds a finite region of $G^-$ since $G$ is theta-connected. 
For $1\le i\le 5$, let $r_i$ be the second region of $G^-$
incident with the edge $w_iw_{i+1}$, where $w_6$ means $w_1$. 
Now $r_4\ne r_5$, since $h_5$ is not a cut-edge of $G^-$, and so
$r_1,r_2,r_3$ are the only regions of $G^-$ incident with two of $h_1\l h_4$.
If some finite region of $G^-$ is incident with three of $h_1\l h_4$, then two of $r_1,r_2,r_3$ are equal and finite,
contradicting the theta-connectivity of $G$. By (5), $r_1,r_3$ are finite regions of $G^-$, and there is an edge of $G^-$ 
incident with $r_1,r_3$, again contrary to the theta-connectivity of $G$. This proves (6).

\bigskip

Let $z_1\l z_8$
be as in \ref{properties},
and for $1\le i\le 8$, let $Z_i$ be the path of $Z$ between $z_i, z_{i+1}$
containing no other vertex in $\{z_1\l z_8\}$ (where $z_9$ means $z_1$).
\\
\\
(7) {\em $W$ contains two of $g_1\l g_4$.}
\\
\\
For suppose $W$ only contains one. From the symmetry we may assume that $g_1\in E(W)$. Consequently $W\setminus g_1$ 
is a four-edge path of $G^-$ between $z_1,z_3$. Since $z_1,z_3$ have degree two in $G^-$, the first and last edges of this path 
belong to the cycle $Z$. If the edge of $W\setminus g_1$ incident with $z_1$ belongs to $Z_8$, then there is a
three-edge path in $G^-$ between $Z_8\setminus z_1$ and $z_3$, contrary to the theta-connectivity of $G$; so the first edge
of $W\setminus g_1$  belong to $Z_1$ and similarly the last edge belongs to $Z_2$, respectively. Since $G$ is theta-connected,
it follows that the middle vertex of the path $W\setminus g_1$ also belongs to $Z_1\cup Z_2$; and so
$W\setminus g_1=Z_1\cup Z_2$. In particular, $z_2\in V(W)$, and so $g_2\in \delta_G(V(W))$. Since every edge
in $F$ belongs to $E(G^-)$, it follows that $g_2=h_5$, and so $z_2=w_5$. Now since $G$ is theta-connected,
there are five edge-disjoint paths of $G$ between $V(W)$ and $Z_5\cup Z_6\cup Z_7$. One of them uses $g_2$, but the other four are
paths of $G^-$, and start at distinct vertices of $W$. Their first edges are the four edges in $F$.
Let these paths be $P_1\l P_4$ say, numbered so $P_i$ has first vertex $w_i$ and first edge $h_i$.
It follows that no region of $G^-$ is incident with three of $h_1\l h_4$, because of the paths
$P_1\l P_4$; and similarly no two regions with a common edge are together incident with all of $h_1\l h_4$, contrary to (5).
This proves (7).

\bigskip

Since $g_1\l g_4$ are pairwise vertex-disjoint, and $W$ has length five, it follows that $W$ contains exactly two
of $g_1\l g_4$. Let $W$ contain $g_1$ and $g_i$ say. The other three edges of $W$ are edges of $G^-$ incident with vertices 
in $\{z_1\l z_8\}$, and so all three belong to $E(Z)$. Consequently $i=2$, and one of $Z_1,Z_3$ has length one, and the
other has length two, and from the symmetry we may assume that $Z_1$ has length one and $Z_3$ length two. From the 
theta-connectivity of $G$, there are five edge-disjoint paths of $G^-$ from $V(W)$ to $Z_5\cup Z_6\cup Z_7$, and their first edges
are the five edges in $\delta_G(V(W))$. Again this contradicts (5), and completes the proof of \ref{mainthm}.~\bbox

\includepdf[pages={1-11}, pagecommand=\thispagestyle{fancy},
scale=0.9]{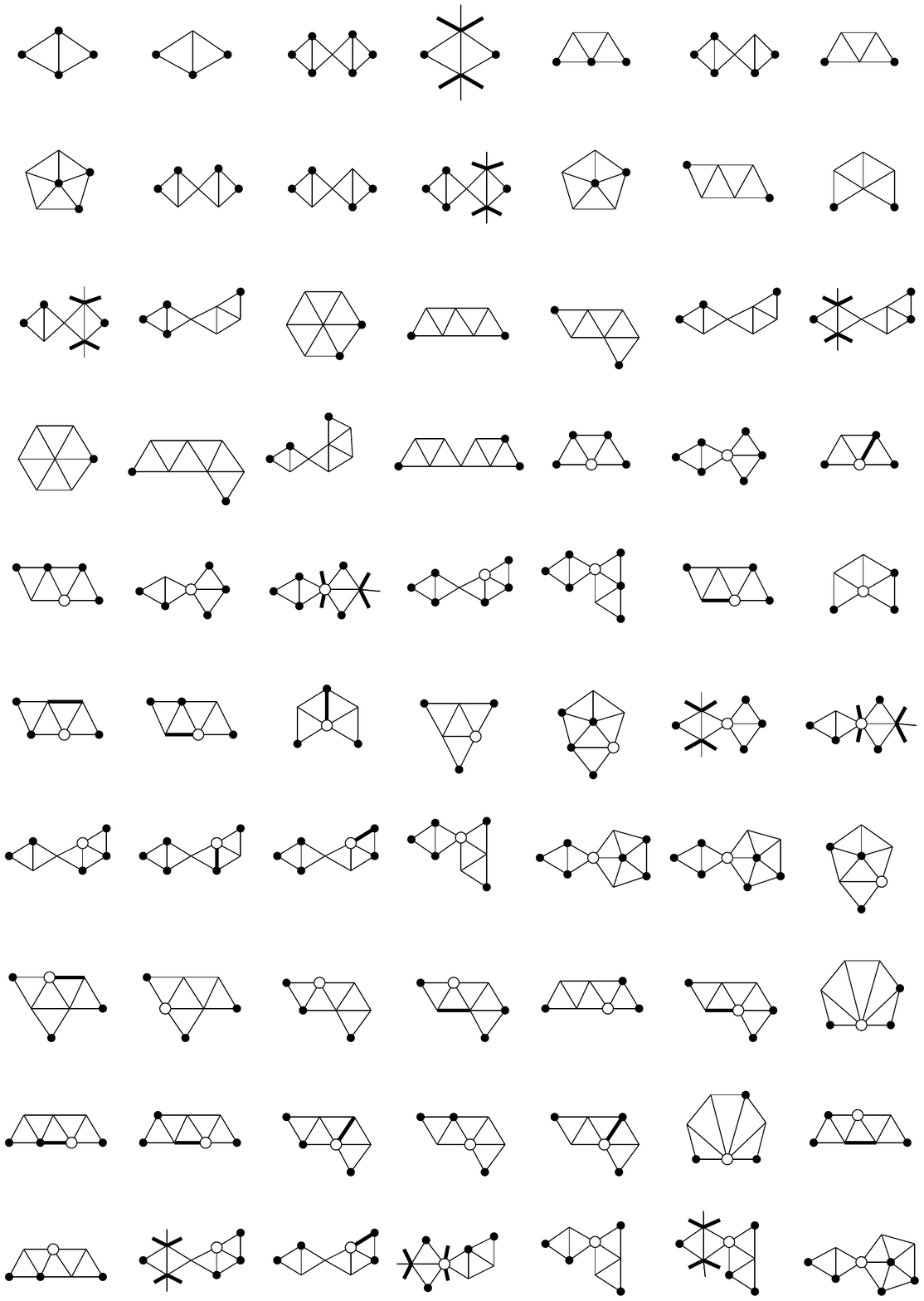}
\end{document}